\documentclass[12pt,a4paper,reqno]{amsart}
\usepackage{amsmath,amsthm}
\usepackage{amsfonts}
\usepackage{amssymb}
\allowdisplaybreaks

\newcommand{\be}{\begin{equation}}
\newcommand{\ef}{\end{equation}}
\chardef\bslash=`\\ % p. 424, TeXbook
\newtheorem*{thm*}{Theorem}

\theoremstyle{definition}

\newtheorem*{remark*}{Remarks}
\newtheorem*{defn*}{Definition}

\theoremstyle{remark}
\newcommand{\G}{\Gamma}
\newcommand{\wt}{\widetilde}
\newcommand{\wh}{\widehat}
\newcommand{\fc}{\frac}

\newcommand{\iy}{\infty}
 \renewcommand{\sectionmark}[1]{}

\newcommand{\Be}{Beltrami}

\newcommand{\qc} {quasiconformal}

\newcommand{\ve}{\varepsilon}

 \usepackage{amsfonts}
\newcommand{\field}[1]{\mathbb{#1}}
\newcommand{\g}{\gamma}

\newcommand{\D}{\field{D}}
\newcommand{\om}{\omega}
\newcommand{\z}{\zeta}
\newcommand{\ov}{\overline}
\newcommand{\vp}{\varphi}
\newcommand{\hC}{\wh{\field{C}}}
\newcommand{\C}{\field{C}}

\newcommand{\B}{\mathbf{B}}
\newcommand{\T}{\mathbf{T}}
\newcommand{\Belt}{\operatorname{Belt}}

\newcommand{\Fib}{\operatorname{Fib}}
\newcommand{\Pot}{\operatorname{Pot}}
\newcommand{\Teich}{\operatorname{Teich}}

\newcommand{\af}{\operatorname{af}}

\newcommand{\x} {\mathbf x}
\renewcommand{\a} {\alpha}

\newcommand{\ld}{\lambda}
\newcommand{\kp}{\kappa}

\begin{document}

\title{On Zalcman's and Bieberbach conjectures}
\author{Samuel L. Krushkal}

\begin{abstract}The well-known Zalcman conjecture, which implies the Bieberbach conjecture, states that the coefficients of univalent functions
$f(z) = z + \sum\limits_2^\iy a_n z^n$ on the unit disk satisfy
$|a_n^2 - a_{2n-1}| \le (n-1)^2$ for all $n > 2$,
with equality only for the Koebe function and its rotations.
The conjecture was proved by the author for $n \le 6$ (using geometric arguments related to the Ahlfors-Schwarz lemma) and remains open for $n \ge 7$.

The main theorem of this paper states that these conjectures are equivalent and provides their
simultaneous proof for all $n \ge 3$  combining
 the indicated geometric arguments with a new  author's approach to extremal problems for holomorphic functions based on lifting the rotationally homogeneous coefficient functionals to the Bers fiber space over universal Teichm\"{u}ller space.
\end{abstract}

\date{\today\hskip4mm({ZaBiCos.tex})}

\maketitle

\bigskip

{\small {\textbf {2020 Mathematics Subject Classification:} Primary:
30C50, 30C75, 30F60; Secondary 30C55, 30C62, 31A05, 32Q45, 32L81}

}

\bigskip

\markboth{Samuel L. Krushkal}{On Zalcman's and Bieberbach conjectures}
\pagestyle{headings}

\bigskip\bigskip
\centerline{\bf 1. INTRODUCTORY REMARKS AND STATEMENT OF RESULT}

\bigskip\noindent
{\bf 1.1. Zalcman's conjecture}.
The holomorphic functionals on the classes of univalent functions depending on the Taylor coefficients of these functions play an important role in various geometric and physical applications of complex analysis, for example, in view of their connection with string theory and with a holomorphic extension of the Virasoro algebra.
These coefficients reflect the fundamental intrinsic features of conformal maps. Thus estimating these coefficients was an important problem in geometric function theory already more then hundred years.

There are two general canonical classes $S$ and $\Sigma$ formed by univalent functions
$$
f(z) = z + a_2 z^2 + a_3 z^3 +\dots \quad \text{and} \ \ F(z) = z + b_0 + b_1 z^{-1} + b_2 z^{-2} + \dots
$$
on the disks $\D = \{|z| < 1\}$ and $\D^* = \{z \in \hC = \C \cup \{\iy\}: \ |z| > 1\}$, which play a crucial  role in geometric function theory,  Teichm\"{u}ller space theory and other fields.

There were several classical conjectures about the coefficients of these functions.
They include the Bieberbach conjecture that in the class $S$  the coefficients are estimated by $|a_n| \le n$, as well as several other well-known conjectures that imply the Bieberbach conjecture.
Most of them have been proved by the de Branges theorem, see \cite {DB}, \cite{Ha},
\cite{Kr6}.

\bigskip
In the 1960s, Lawrence Zalcman posed the conjecture that {\it for
any $f \in S$ and all $n \ge 3$,
 \be\label{1}
|a_n^2 - a_{2n-1}| \le (n-1)^2,
\end{equation}
with equality only for the Koebe function}
 \be\label{2}
\kp_\theta(z) = \fc{z}{(1 - e^{i \theta} z)^2} = z +
\sum\limits_2^\iy n e^{- i(n-1) \theta} z^n, \quad - \pi < \theta \le \pi,
\end{equation}
which maps the unit disk onto the complement of the ray
$$
w = -t e^{-i \theta}, \ \ \fc{1}{4} \le t \le \iy.
$$
This remarkable conjecture implies the Bieberbach conjecture. The exceptional case $n = 2$ is simple and somewhat different: the equality $|a_2^2 - a_3| = 1$ is valid for the Koebe function $\kp_\theta(z)$ and the odd function
$\sqrt{\kp_\theta(z^2)} = z/(1 - e^{i\theta} z^2)$ (note that the existence of extremals for
$Z_n(f) = a_n^2 - a_{2n-1}$  follows from compactness of the class $S$ with respect to locally uniform convergence on $\D$).

\bigskip
Zalcman's conjecture has been proved by the author for $n = 3, 4, 5, 6$ in \cite{Kr2}, \cite{Kr3} by applying the geometric arguments related to the Ahlfors-Schwarz lemma and its generalizations;
it remains open for $n \ge 7$.

In the papers \cite{BT}, \cite{Ma1}, this conjecture was proved for certain special subclasses of $S$, and later their results have been generalized by several authors, which considered somewhat generalized Zalcman's functional (see, e.g., \cite{TTV} and the references there).
Eframidis and Vukoti\'{c} proved the asymptotic version of Zalcman's conjecture using the famous Hayman theorem on the asymptotic growth of coefficients of individual functions
$f \in S$ (see, \cite{EV}, \cite{Ha}). Recently Obradovi\'{c} and Tuneski \cite{OT} proved
the conjecture for univalent functions with real coefficients.

This conjecture is now restated (for $n \ge 7$) in the book \cite{HL} collecting the research problems in the function theory.

\bigskip\noindent
{\bf 1.2. Main theorem}. Recently the author established in \cite{Kr7}-\cite{Kr9} a new approach to solving the classical coefficient problems on various classes of holomorphic functions, not necessarily univalent.
This approach is based on lifting the given polynomial (and more general holomorphic) rotationally homogeneous functionals
$$
J(f) = J(a_{m_1}, \dots, a_{m_s})
$$
depending from the distinguished finite set of coefficients $a_{m_j}$ satisfying
$$
2 < a_{m_1} < \dots < a_{m_s} < \iy
$$
onto the Bers fiber space over the universal Teichm\"{u}ller space and applying the deep analytic and geometric features of Teichm\"{u}ller spaces, especially the Bers isomorphism theorem for Teichm\"{u}ller spaces of punctured Riemann surfaces.

One obtains on this way, for example, an alternate and direct proof of the Bieberbach conjecture and even
of more stronger conjectures, estimating the Taylor coefficients of the Schwarzian derivatives, of the Hardy functions, etc.

The arguments applied in \cite{Kr7}-\cite{Kr9} essentially used the assumption that functional $J(f)$ is  strongly rotationally homogeneous, which means homogeneity under the pre and post rotations of functions
 \be\label{3}
f(z) \mapsto f_{\a, \beta}(z) = e^{i \beta} f(e^{i\a} z)
\end{equation}
with different $\a$ and $\beta$ from $[- \pi, \pi)$. However, Zalcman's functional $Z_n(f)$ admits only
a weaker rotational homogeneity (with $\beta = - \a$).

The aim of this paper is to show that combining both indicated approaches (with some needed modification)  one obtains the following theorem. which proves these conjectures simultaneously and establishes their equivalence.

Fix $n \ge 3$ and consider on the class $S$ the coefficient functionals
 \be\label{4}
J_1(f) = \fc{a_n^2 - a_{2n-1}}{(n - 1)^2}, \quad J_2(f) = \fc{a_n}{n}.
\end{equation}
dominated by subharmonic functional
 \be\label{5}
J_{1,2}(f) = \max \left \{\max_S |J_1(f)|, \ \max_S |J_2(f)| \right\}.
\end{equation}

\bigskip\noindent
{\bf Theorem 1}. {\it For any function $f \in S$ and any $n \ge 3$, we have the sharp bound
$|J_{1,2}(f)| \le 1$, with equalities
  \be\label{6}
\max_S |J_1(f)| = \max_S |J_2(f)|= 1
\end{equation}
only for the Koebe function (2) and its rotations $e^{- i\a} \kp_\theta(e^{i \a} z)$.
}

\bigskip\noindent
The proof of this theorem  leads to the bound $J_{1,2}(f) \le 1$ for all $f \in S$,
with the equality $J_{1,2}(f_0) = 1$ for any extremal function maximizing $J_{1,2}$. This yields by (5) that
either $\max_S |J_1(f_0)| = 1$ and $\max_S |J_2(f_0) \le 1$,  or $\max_S |J_1(f_0) \le 1$ and $\max_S |J_2(f_0) = 1$, or both equalities (6). In either case the extremality is possible only for $f_0 = \kp_\theta$ (and then both these maxima must be equal to $1$).

\bigskip\noindent
{\bf 1.3. On extremality of Koebe's function}. 
The spaces $\mathcal T_m$ possess an important property
inherited from the classes  $\Sigma$ and $\wh \Sigma(1)$ inherited from the class $\Sigma$
and preserved under connection (1). This is a circular symmetry of maps $F_{\tau,\sigma}$,
since for any $F \in \Sigma$, all functions $e^{- i \a} F(e^{i \a}z)$ also belong to this
class. Theorem 1 intrinsically relates to the general result on extremality of Koebe's function established in \cite{Kr10} concerning the arbitrary rotationally homogeneous polynomial functionals
$J(f) = J(a_{m_1}, \dots, a_{m_s})$ on the class $S$ (with $2 < a_{m_1} < \dots < a_{m_s} < \iy$).

The relations  between the coefficients of $f(z) = z + a_2 z^2 + \dots \in S$ and of their inversions  $F(z) = 1/f(1/z) = z + b_0 + b_1 z^{-1} + \dots$, which belong to the class $\Sigma$ of univalent $\hC$-holomorphic functions on the disk
$\D^* = \{z \in \hC = \C \cup \{\iy\}: \ |z| > 1\}$ with such expansions,
transform $J(f)$ to a functional $\wt J(F)$ on $\Sigma$, which is called to be associated with $J$ (see below the relations (20)). Denote by $\Sigma_{\af}$ the subcollection in $\Sigma$ formed by functions
$$
F_{b_0,b_1;t}(z) = z + b_0 t + b_1 t^2 z^{-1}
$$
with
$$
|b_0| \le 2, \ \ |b_1| \le 1, \ \ |t| \le 1.
$$
These functions have the affine extensions $\wh F_{b_0,b_1;t}(z) = z + b_0 t + b_1 t^2 \ov z$ onto
the unit disk $\D = \{|z| < 1\}$ (with constant dilatations $b_1t^2$) .

The indicated result states:

\bigskip\noindent
{\bf Theorem 2}. {\it The Koebe function $\kp_\theta(z)$ is (a unique) extremal of a rotationally homogeneous coefficient functional $J(f)$ if and only if its associated functional $\wt J(F)$ satisfies
$$
\max_\Sigma |\wt J(F)| = \max_{f \in S} |\wt J(F_f)| = \sup_{\Sigma_{\af}} |\wt J(F_{b_0,b_1;t})|.
$$

In other words, the maximal value of the associate to $J$ functional $\wt J$ on $\Sigma$ must be  attained  on the distinguished subset $\Sigma_{\af}$ of functions $F \in \Sigma$ admitting the affine
extensions to the disk $\D$.  }

\bigskip
This theorem shows that actually the collection of coefficient functionals maximized by Koebe's function
is rather sparse. Giving rise to possible generalizations of the Bieberbach and Zalcman conjectures, it simultaneously describes the admissible extent of this direction.

Note that this theorem involves only the weak rotational homogeneity of $J(f)$.

\bigskip\bigskip
\centerline{\bf 2. DIGRESSION TO TEICHM\"{U}LLER SPACES}

\bigskip\noindent
{\bf 2.1}.
First, we briefly recall some needed results from Teichm\"{u}ller space theory involved in the proof of  Theorem 1; the details can be found, for example, in \cite{Be2}, \cite{GL}, \cite{Le}.

This theory is intrinsically connected with univalent functions having quasiconformal extensions onto $\hC$, which involves three normalization conditions to have uniqueness, compactness, holomorphic dependence on parameters, etc. It is technically more convenient to deal with functions univalent on the disk $\D^*$.

The (dense) subcollections of functions $f \in S$ and $F \in \Sigma$ admitting quasiconformal extensions across the unit circle $\mathbb S^1 = \{|z| = 1\}$ (and hence onto $\hC$) will be denoted by $S_Q$ and $\Sigma_Q$, respectively.

\bigskip\noindent
{\bf 2.2}.
The {\bf universal Teichm\"{u}ller space} $\T = \Teich (\D)$ is the space of quasisymmetric homeomorphisms of the unit circle $\mathbb S^1$ factorized by M\"{o}bius maps;  all Teichm\"{u}ller spaces have their isometric
copies in $\T$.

The canonical complex Banach structure on $\T$ is defined by factorization of the ball of
the {\bf Beltrami coefficients} (or complex dilatations)
$$
\Belt(\D)_1 = \{\mu \in L_\iy(\C): \ \mu|\D^* = 0, \ \|\mu\| < 1\},
$$
letting $\mu_1, \mu_2 \in \Belt(\D)_1$ be equivalent if the
corresponding \qc \ maps $w^{\mu_1}, w^{\mu_2}$ (solutions to the
Beltrami equation $\partial_{\ov{z}} w = \mu \partial_z w$ with $\mu
= \mu_1, \mu_2$) coincide on the unit circle $\mathbb S^1 = \partial \D^*$
(hence, on $\ov{\D^*}$). Such $\mu$ and the corresponding maps
$w^\mu$ are called $\T$-{\it equivalent}. The equivalence classes
$[w^\mu]_\T$ are in one-to-one correspondence with the {\bf Schwarzian derivatives}
$$
S_w(z) = \left(\frac{w^{\prime\prime}(z)}{w^\prime(z)}\right)^\prime
- \frac{1}{2} \left(\frac{w^{\prime\prime}(z)}{w^\prime(z)}\right)^2
\quad (w = w^\mu(z), \ \ z \in \D^*).
$$

The solutions of the Beltrami equation $\partial_{\ov{z}} w = \mu \partial_z w$ with $\mu \in \Belt(\D)_1$ and of the Schwarz equation $S_w(z) = \vp(z)$ with $\vp \in \B(\D^*)$ are conformal
on the disk $\D^*$ and are defined up to a M\"{o}bius (fractional-linear) transformation of $\hC$.
We shall use the hydrodynamical normalization on the infinite point $z = \iy$, i.e., let
$$
w^\mu(z) = z + b_0 + b_1 z^{-1} + b_2 z^{-2} + \dots
$$
and some additional normalization, for example $w^\mu(1) = 1$.
The chain rule
$$
S_{F_1 \circ F}(z) = (S_{F_1} \circ F) F^\prime(z)^2 + S_F(z).
$$
yields for the M\"{o}bius map $w = \g(z)$ the equality
$$
S_{F_1 \circ \g}(z) = (S_{f_1} \circ \g) \g^\prime(z)^2,
\quad S_{\g \circ F}(z) = S_F(z).
$$
Hence, each $S_w(z)$ can be regarded as a quadratic differential $\vp = S_w(z) dz^2$ on $\D^*$.
The solution $w(z)$ of the Schwarzian equation $S_w(z) = \vp(z)$ with a given holomorphic $\vp$ is defined up to a M\"{o}bius transformation of $\hC$.

For every locally univalent function $w(z)$ on a simply connected hyperbolic domain $D \subset \hC$, its Schwarzian derivative belongs to the complex Banach space $\B(D)$ of hyperbolically bounded holomorphic functions on $D$ with the norm
$$
\|\vp\|_{\B(D)} = \sup_D \ld_D^{-2}(z) |\vp(z)|,
$$
where $\ld_D(z) |dz|$ is the hyperbolic metric on $D$ of Gaussian curvature $- 4$; hence, $\vp(z) = O(z^{-4})$ as $z \to \iy$ if $\iy
\in D$. In particular,
$$
\ld_\D(z) = 1/(1 - |z|^2), \quad \ld_{\D^*}(z) = 1/(|z|^2 -1).
$$
The space $\B(D)$ is dual to the Bergman space $A_1(D)$, a subspace
of $L_1(D)$ formed by integrable holomorphic functions (quadratic differentials
$\vp(z) dz^2$) on $D$.

The Schwarzians $S_{w^\mu}(z)$ with $\mu \in \Belt(\D)_1$ range over a bounded domain in the space $\B = \B(\D^*)$, which models the space $\T$.
It is located in the ball $\{\|\vp\|_\B < 6\}$ and contains the ball $\{\|\vp\|_\B < 2\}$.
In this model, the Teichm\"{u}ller spaces of all hyperbolic Riemann surfaces are contained in $\T$ as its complex submanifolds.

The factorizing projection $\phi_\T(\mu) = S_{w^\mu}: \ \Belt(\D)_1 \to \T$
is a holomorphic map from $L_\iy(\D)$ to $\B$. This map is a split
submersion, which means that $\phi_\T$ has local holomorphic sections (see, e.g., \cite{Be1}, \cite{EKK}, \cite{GL}).

Both equations $S_w = \vp$ and $\partial_{\ov z} w = \mu \partial_z w$ (on $\D^*$ and $\D$, respectively) determine their solutions up to a M\"{o}bius transformation of $\hC$.

The following important lemma from \cite{Kr8} provides a somewhat different normalization of quasiconformally extendable functions, which also insures the needed uniqueness of solutions, their holomorphic dependence of complex parameters, etc.

\bigskip\noindent
{\bf Lemma 1}. {\it For any Beltrami coefficient $\mu \in \Belt(\D^*)_1$ and any $\theta_0 \in [0, 2 \pi]$, there exists a point $z_0 = e^{i \a}$ located on $\mathbb S^1$ so that
$|e^{i \theta_0} - e^{i \a}| < 1$ and such that for any $\theta$ satisfying
$|e^{i \theta} - e^{i \a}| < 1$ the equation
$\partial_{\ov z} w =  \mu(z) \partial_z w$
has a unique homeomorphic solution $w = w^\mu(z)$, which is holomorphic on the unit disk $\D$
and satisfies
 \be\label{7}
w(0) = 0, \quad w^\prime(0) = e^{i \theta}, \quad w(z_0) = z_0.
\end{equation}
Hence, $w^\mu(z)$ is conformal and does not have a pole in $\D$ \ (so
$w^\mu(z_{*}) = \iy$ at some point $z_{*}$ with $|z_{*}| \ge 1$).  }

\bigskip
In particular, this lemma allows one to define the Teichm\"{u}ller spaces using the quasiconformally extendible  univalent functions $w(z)$ in the unit disk $\D$ normalized by
 \be\label{8}
w(0) = 0, \quad w^\prime(0) = 1, \quad w(1) = 1
\end{equation}
and with more general normalization
$$
w(0) = 0, \quad w^\prime(0) = e^{i \theta}, \quad w(1) = 1.
$$
All such functions are holomorphic in the disk $\D$.

As a simple consequence of this lemma, we have, passing to rotations $w(z) \mapsto e^{i \a} w(e^{-i \a} z)$, the following

\bigskip\noindent
{\bf Lemma 2}. {\it For any $\mu \in \Belt(\D^*)_1$, any given value $\theta \in (- \pi, \pi]$ and
any point $z_0$ of the unit circle $\mathbf S^1 = \{|z| = 1\}$,
there exists a unique homeomorphic solution $w = w^\mu(z)$ of the equation
$\partial_{\ov z} w = \mu(z) \partial_z w$ on $\hC$ such that
 \be\label{9}
w(0) = 0, \quad w^\prime(0) = 1, \quad w(z_0) = e^{i \theta}.
\end{equation}
This solution is holomorphic on the unit disk $\D$, and hence, $w(z_{*}) = \iy$ at some point $z_{*}$ with}  $|z_{*}| \ge 1$.

\bigskip
Note that for $\mu(z) = 0$ (almost everywhere on $\D^*$) the corresponding solution $w^\mu(z)$ with
$w(0) = 0, \quad w^\prime(0) = e^{i \theta}, \quad w(1) = 1$
is the elliptic M\'{o}bius map
$$
w = \fc{e^{- i \theta} z}{(e^{- i \theta} - 1) z + 1}
$$
(with the fixed points $0$ and $1$); it equals the identity map if $\theta = 0$.

\bigskip\noindent
{\bf 2.3}. The points of {\bf Teichm\"{u}ller space $\T_1 = \Teich(\D_{*})$ of the punctured disk} $\D_{*} = \D \setminus \{0\}$ are the classes $[\mu]_{\T_1}$ of $\T_1$-{\bf equivalent} Beltrami coefficients $\mu \in \Belt(\D)_1$, which means that the corresponding quasiconformal  automorphisms $w^\mu$ of the unit disk coincide on both boundary components of $\D_{*}$ (the unit circle $\mathbb S^1 = \{|z| =1\}$ and the puncture $z = 0$) and are homotopic on $\D \setminus \{0\}$. This space can be endowed with a canonical complex structure of a complex Banach manifold and embedded into $\T$ using uniformization.

Namely, the disk $\D_{*}$ is conformally equivalent to the factor $\D/\G$, where $\G$ is a cyclic parabolic Fuchsian group acting discontinuously on $\D$ and $\D^*$. The functions $\mu \in L_\iy(\D)$ are lifted to $\D$ as the \Be \ $(-1, 1)$-measurable forms  $\wt \mu d\ov{z}/dz$ in $\D$ with respect to $\G$, i.e., via
$(\wt \mu \circ \g) \ov{\g^\prime}/\g^\prime = \wt \mu, \ \g \in \G$, forming the Banach space $L_\iy(\D, \G)$.

We extend these $\wt \mu$ by zero to $\D^*$ and consider the unit ball $\Belt(\D, \G)_1$ of $L_\iy(\D, \G)$. Then the corresponding Schwarzians $S_{w^{\wt \mu}|\D^*}$ belong to $\T$. Moreover, $\T_1$
is canonically isomorphic to the subspace $\T(\G) = \T \cap \B(\G)$, where $\B(\G)$ consists of elements $\vp \in \B$ satisfying $(\vp \circ \g) (\g^\prime)^2 = \vp$ in $\D^*$ for all $\g \in \G$.

Due to the {\bf Bers isomorphism theorem}, {\it the space $\T_1$ is biholomorphically isomorphic to the Bers fiber space
$$
\mathcal F(\T) = \{(\phi_\T(\mu), z) \in \T \times \C: \ \mu \in
\Belt(\D)_1, \ z \in w^\mu(\D)\}
$$
over the universal space $\T$ with holomorphic projection} $\pi(\psi, z) = \psi$ (see \cite{Be2}).

This fiber space is a bounded hyperbolic domain in $\B \times \C$ and represents the collection of domains $D_\mu = w^\mu(\D)$ as a holomorphic family over the space $\T$. For every $z \in \D$,  its
orbit $w^\mu(z)$ in $\T_1$ is a holomorphic curve over $\T$.

The indicated isomorphism between $\T_1$ and $\mathcal F(\T)$ is induced by the inclusion map \linebreak
$j: \ \D_{*} \hookrightarrow \D$ forgetting the puncture at the origin via
 \be\label{10}
\mu \mapsto (S_{w^{\mu_1}}, w^{\mu_1}(0)) \quad \text{with} \ \
\mu_1 = j_{*} \mu := (\mu \circ j_0) \ov{j_0^\prime}/j_0^\prime,
\end{equation}
where $j_0$ is the lift of $j$ to $\D$.

The Bers theorem is valid for Teichm\"{u}ller spaces $\T(X_0 \setminus \{x_0\})$ of all punctured hyperbolic Riemann surfaces $X_0 \setminus \{x_0\}$ and implies that $\T(X_0 \setminus \{x_0\})$ is
biholomorphically isomorphic to the Bers fiber space $ \Fib(\T(X_0))$ over $\T(X_0)$.

\bigskip\noindent
{\bf 2.4}. The spaces $\T$ and $\T_1$ can be weakly (in the topology generated by the spherical metric on $\hC$) approximate by finite dimensional Teichm\"{u}ller spaces $\T(0, n)$ of punctured spheres (Riemann surfaces of genus zero)
$$
X_{\mathbf z} = \hC \setminus \{0, 1, z_1 \dots, z_{n-3}, \iy\}
$$
defined by ordered $n$-tuples $\mathbf z = (0, 1, z_1, \dots, z_{n-3}, \iy), \ n > 4$ with distinct
$z_j \in \C \setminus \{0, 1\}$ (the details see, e.g., in \cite{Kr7}).

Fix a collection $\mathbf z^0 = (0, 1, z_1^0, \dots, z_{n-3}^0, \iy)$ with $ z_j^0 \in \mathbb S^1$
defining the base point $X_{\mathbf z^0}$ of the space $\T(0, n) = \T(X_{\mathbf z^0})$. Its points are the equivalence classes $[\mu]$ of Beltrami coefficients from the ball
$\Belt(\C)_1 = \{\mu \in L_\iy(\C): \ \|\mu\|_\iy < 1\}$ under the relation: $\mu_1 \sim \mu_2$, if the corresponding quasiconformal  homeomorphisms $w^{\mu_1}, w^{\mu_2}: \ X_{\mathbf a^0} \to X_{\mathbf a}$  are homotopic on $X_{\mathbf a^0}$ (and hence coincide in the points $0, 1, z_1^0, \dots, z_{n-3}^0, \iy$).
This models $\T(0, n)$ as the quotient space $\T(0, n) = \Belt(\C)_1/\sim $ with complex Banach structure of dimension $n - 3$ inherited from the ball $\Belt(\C)_1$.

Another canonical model of the space $\T(0, n) = \T(X_{\mathbf z^0})$ is obtained again using the uniformization. The surface $X_{\mathbf z^0}$ is conformally equivalent to the quotient space $U/\G_0$, where
$\G_0$ is a torsion free Fuchsian group of the first kind acting discontinuously on $\D \cup \D^*$. The functions $\mu \in L_\iy(X_{\mathbf z^0})$ are lifted to $\D$ as the Beltrami $(-1, 1)$-measurable forms  $\wt \mu d\ov{z}/dz$ in $\D$ with respect to $\G_0$ which satisfy
$$
(\wt \mu \circ \g) \ov{\g^\prime}/\g^\prime = \wt \mu, \ \g \in \G_0,
$$
and form the Banach space $L_\iy(\D, \G_0)$.
After extending these $\wt \mu$ by zero to $\D^*$, the  Schwarzians $S_{w^{\wt \mu}|\D^*}$ for
$\|\wt \mu\|_\iy < 1$  belong to $\T$ and form its subspace regarded as the {\bf Teichm\"{u}ller space
$\T(\G_0)$ of the group $\G_0$}. It is canonically isomorphic to the space $\T(X_{\mathbf z^0})$, and moreover,
$$
\T(\G_0) = \T \cap \B(\G_0),
$$
where $\B(\G_0)$ is an $(n - 3)$-dimensional subspace of $\B$ which consists of elements $\vp \in \B$ satisfying $(\vp \circ \g) (\g^\prime)^2 = \vp$ for all $\g \in \G_0$ (holomorphic $\G_0$-automorphic forms of degree $- 4$); see, e.g. \cite{Le}. This space has has the same elements as the space $A_1(\D^*,\G_0)$ of integrable holomorphic forms of degree $- 4$.

This leads to the representation of the space $\T(X_{\mathbf z^0})$ as a bounded domain in the complex Euclidean space $\C^{n-3}$.

Any Teichm\"{u}ller space is a complete metric space with intrinsic Teichm\"{u}ller metric defined by quasiconformal maps. By the Royden-Gardiner theorem, this metric equals the hyperbolic  Kobayashi metric
determined by the complex structure (see, e.g., \cite{EKK}, \cite{GL}, \cite{Ro1}).

\bigskip\bigskip
\centerline{\bf 3. PRELIMINARIES AND UNDERLYING GEOMETRIC RESULTS}

\bigskip\noindent
{\bf 3.1. Holomorphic homotopy of a univalent function}.
For each $F(z) = z + b_0 + b_1 z^{-1} + \dots \in \Sigma$
we define a complex homotopy
 \be\label{11}
F_t (z) = t F \left( \fc{z}{t} \right) =  z + b_0 t + b_1 t^2 z^{-1} +
b_2 t^3 z^{-2} + ...: \ \D^* \times \D \to \hC
\end{equation}
of this function to the identity map. Then
$$
S_{F_t}(z) = t^{-2} S_F(t^{-1} z)
$$
and, moreover, the map $h_F: \  t \mapsto S_{F_t}$ is holomorphic as a function $\D \to \B$ (see, e.g. \cite{Kr2}). It determines the homotopy disk $\D(F) = \{F_t\}$, which is holomorphic at the noncritical points of $h_F$.
These disks foliate the set $\Sigma_Q$.

The corresponding homotopy of functions from $S$ is given by
$f(z, t) = t^{-1} f(t z) = z + a_2 t + \dots$; so $Z_n(f(\cdot, t)) = t^{2n-2} Z_n(f)$.

The dilatation of the homotopy functions is estimated by the following lemma, which is a special case
of dynamical $r^2$-property of domains and univalent functions (see \cite{Kr2}, \cite{KK2}).

\bigskip\noindent
{\bf Lemma 3}. \cite{Kr2} {\it (a) Each homotopy map $F_t$ admits $k$-quasiconformal extension to the whole sphere $\hC = \C \cup \{\iy\}$ with $k \le |t|^2$.
The bound $k(F_t) \le |t|^2$ is sharp and occurs only for the maps
$$
F_{b_0,b_1;1} (z) = z + b_0 + b_1 z^{-1}, \quad |b_1| = 1,
$$
whose homotopy maps
 \be\label{12}
F_{b_0,b_1;t}(z) = z + b_0 t + b_1 t^2 z^{-1}
\end{equation}
have the affine extensions $\wh F_{b_0,b_1;t}(z) = z + b_0 t + b_1 t^2 \ov z$ onto $\D$.

(b) If $F(z) = z + b_0 + b_p z^{-p} + b_{p+1} z^{-(p+1)} + \dots \ (b_p \ne 0)$
for some integer $p > 1$, then the minimal dilatation of extensions is estimated by $k(F_t) \le |t|^{p+1}$; this bound also is sharp.
}

\bigskip
This lemma is a special case of the dynamical $r^2$-property of domains and univalent functions
investigated in  \cite{Kr2}, \cite{KK2}).

If $F \in \Sigma_Q$ admits a $k$-quasiconformal extension, then the best bound for the dilatation of its homotopy $F_t$ much stronger. This bound will not be applied here.

Note also that, due to Strebel's frame mapping condition \cite{St}, the extremal extensions $\wh F_t$ of any homotopy functions $F_t$ with $|t| <1$ is of Teichm\"{u}ller type, i.e., with the Beltrami coefficient of the form
$$
\mu_{\wh F_t}(z) = \tau(t) |\psi(z)|/\psi(z),
$$
where $\psi$ is a holomorphic function from $L_1(\D)$ (and unique).

\bigskip\noindent
{\bf 3.2}. It suffices for our goals to consider the functions $f \in S$ with
$$
b_1 = a_2^2 - a_3 \ne 0;
$$
this assumption is equivalent to
$S_f(0) = \lim\limits_{z\to \iy} z^4 S_{F_f}(z) \ne 0$. Such functions form a dense subset of $S$, and
their Schwarzians form a dense subset of the space $\T$.

We divide every homotopy function $F_t$ of $ F = F_f$ into two parts
$$
F_t(z) = z + b_0 t + b_1 t^2 z^{-1} + b_2 t^3 z^{-2} + \dots
= F_{b_0,b_1;t}(z) + h(z,t),
$$
where $F_{b_0,b_1;t}$ is the map (12) with $b_0, \ b_1$ coming from $F$.
For sufficiently small $|t|$, the remainder $h$ is estimated by
$h(z,t) = O(t^3)$ uniformly in $z$ for all $|z| \ge 1$.

Then the Schwarzian derivatives of $F_t$ and $F_{b_0,b_1;t}$
are related by
$$
S_{F_t}(z) = S_{F_{b_0,b_1;t}}(z) + \om(z, t),
$$
where the remainder $\om$ is uniquely determined by the chain rule
$$
S_{w_1\circ w}(z) = (S_{w_1} \circ w) (w^\prime)^2(z) + S_w(z),
$$
and is estimated in the norm of $\B$ by
and is estimated in the norm of $\B$ by
 \be\label{13}
\|\om(\cdot,t)\|_\B = O(t^3), \quad t \to 0;
\end{equation}
this estimate is uniform for $|t| < t_0$ (cf., e.g. \cite{Be1}, \cite{Kr1}).

All functions $F_{b_0,b_1;t}$ with
\be\label{14}
|b_0| \le 2, \ \ |b_1| \le 1, \ \ |t| \le 1.
\end{equation}
are univalent on the disk $\D^*$ (but can vanish there) and, if $|t| < 1$, have the affine  extensions
onto $\D$. For such functions, their homotopy disk $\D(F) = \{F_t\}$ coincides with the extremal disk
$\D(\psi) = \{t \mu_0: \ t \in \D\} \subset \Belt(\D)_1$; hence, the action of functional $\wt Z_n$ on extermal disks of functions $F_{b_0,b_1;t}$ is rotationally symmetric with respect to $t \in \D$.

\bigskip
We call the values $b_0$ and $b_1$ \textbf{admissible} if they are the initial coefficients of some function from $\Sigma_Q$ (these values satisfy (14)). The collection of all such $F_t$  with $|t| < 1$ will be denoted by $\Sigma_{\af}$.

\bigskip\noindent
{\bf 3.3. Two generalizations of the Gaussian curvature}. The next geometric results of this section
are needed for the Zalcman's part of Theorem 1. It relies on the properties of subharmonic conformal metrics $\ld(z) |dz|$ on the disk (with $\ld(z) \ge 0$), whose curvature is at most $- 4$ in the generalized sense introduced by Ahlfors and Royden (see \cite{Ah}, \cite{He}, \cite{Ro1}).

Recall that the Gaussian curvature of a $C^2$-smooth metric $\ld > 0$ is defined by
$$
\kappa_\ld = - \fc{\D \log \ld}{\ld^2},
$$
where $\D$ means the Laplacian $4 \partial \ov{\partial}$.

A metric $\ld(z) |dz|$ in a domain $G$ on $\C$ (or on a Riemann surface) has curvature less than or equal to $K$ \textbf{in the supporting sense} if for each $K^\prime > K$
and each $z_0$ with $\ld(z_0) > 0$, there is a $C^2$-smooth \textbf{supporting metric} $\wh \ld$ for $\ld$ at $z_0$ (i.e., such that $\wh \ld(z_0) = \ld(z_0)$ and $\wh \ld(z) \le \ld(z)$ in a neighborhood of $z_0$) with $\kappa_{\wh \ld}(z_0) \le K^\prime$, or equivalently,
 \be\label{15}
\D \log \ld \ge - K \ld^2.
\end{equation}

A metric $\ld$ has curvature at most $K$ \textbf{in the potential sense} at $z_0$ if there is a disk
$U$ about $z_0$ in which the function
$$
\log \ld + K \Pot_U(\ld^2),
$$
where $\Pot_U$ denotes the logarithmic potential
$$
\Pot_U h = \fc{1}{2 \pi} \iint\limits_U h(\z) \log |\z - z| d \xi d \eta \quad (\z = \xi + i \eta),
$$
is subharmonic. Since $\D \Pot_U h = h$ (in the sense of distributions), one can replace $U$ by any open subset $V \subset U$, because the function $\Pot_U(\ld^2) - \Pot_V(\ld^2)$ is harmonic on $U$.
The inequality (15) holds for the generic subharmonic metrics also in the sense of distributions.
Note also that the condition of having curvature at most $- K$ in the potential sense
is invariant under conformal maps.

\bigskip\noindent
{\bf Lemma 4}. \cite{Ro1} {\it If a conformal metric has curvature at most $K$
in the supporting sense, then it has curvature at most $K$ in the potential sense.  }

\bigskip
The following lemma concerns the {\bf circularly symmetric} metrics which are the functions of $r = |z|$.

\bigskip\noindent
{\bf Lemma 5}. \cite{Kr6}, \cite{Ro1} {\it Let $\ld(|t|) d |t|$ be a circularly symmetric subharmonic
metric on $\D$ such that
\be\label{16}
\ld(r) = m c r^{m-1} + O(r^{m}) \quad \text{as} \ \ r \to 0 \ \
\text{with} \ \ 0 < c \le 1 \ \ (m = 1, 2, \dots)
\end{equation}
and this metric has curvature at most $- 4$ in the potential sense. Then  }
\be\label{17}
\ld(r) \ge \fc{m c r^{m-1}}{1 - c^2 r^{2m}}.
\end{equation}

\bigskip
In the case $m = 1$ (with $c \ne 0$) applied in \cite{Ro1}, the estimate (17) is reduced to $\ld(r) \ge c/(1 - c^2 r^2)$
with $c = \d(0)$,and right-hand side defines a supporting conformal metric for $\ld$
at the origin with constant Gaussian curvature $- 4$ on the whole disk $\D$.

The dominant of metrics subject to (16) is given by the following

\bigskip\noindent
{\bf Lemma 6}. {\it Let $\ld(t) |dt|$ be a continuous conformal metric on the disk $\D$ with growth (17) near the origin and having the curvature $- 4$ in the supporting sense at its noncritical points. Then}
$$
\ld(t) \le \ld_m(t) = \fc{m |t|^{m-1}}{1 - 2 t^{2m}}  \quad \text{for all} \ \ t \in \D.
$$

\noindent
{\bf 3.3. Restoration of functional $Z_n(f)$ by its infinitesimal form}.
We also mention two important features of Zalcman's functional. Pass to the normalized functional $Z_n^0 = Z_n/\max_S |Z_n(f)|$ mapping $S$ onto the unit disk, and consider its action on infinite holomorphic families
$\mathcal F_S = \{f_t(z) = f(z, t)\} \subset S$ and on $\Sigma_{\af}$
with $F_(z, \cdot) = 1/f(1/z, \cdot), \ t \in \D$.

Using the relations between the coefficients $a_n(t)$ of $f(z, t)$ and the corresponding coefficients     $b_j(t)$ of $F_f(z, t)$, we represent $Z_n^0$ as a polynomial functional on $\Sigma$,
$$
Z_n^0 (f) = \wt Z_n^0(F_f) = \wt Z_n^0(b_0, b_1, \dots, \ , b_{2n-3}).
$$
The given holomorphic families determine the sequences of holomorphic maps
$$
g_m(t) = \wt Z_n(F_m(\cdot,t)): \ \D \to \C, \quad m = 1, 2, \dots \ (F_m \in \Sigma_{\af}),
$$
and their upper envelope
$$
\wh g = \sup_m |g_m(t)|: \ \D \to \D
$$
followed by upper semicontinuous regularization
$\wh g(t) =  \limsup\limits_{t^\prime \to t} g_m(t^\prime)$ presents a logarithmically subharmonic
function on the unit disk.

The maps $g_m$ pull back the hyperbolic metric of this disk $\ld_\D(z)|dz|$ generating on $\D$ the  logarithmically subharmonic  conformal metrics
$ds = \ld_{g_m}(t) |dz|$ with
$$
\ld_{g_m}(t) = g_m^* \ld_\D(t) = \fc{|g_m^\prime(t)|}{1 - |g_m(t)|^2}
$$
of Gaussian curvature $- 4$ at noncritical points of $g_m$. Passing to the upper envelope
 \be\label{18}
\ld_{Z_n^0}(z) = \sup_m \ld_{g_m}(z)
\end{equation}
and its upper semicontinuous regularization, one obtains a logarithmically subharmonic metric on $\D$,
whose curvature is less than or equal to $- 4$ in both supporting and potential senses (cf. \cite{Kr6}).

We shall apply the results on the curvatures indicated above to the values of $\wt Z_n$ on the homotopy disks; thus the derivatives must be understand as distributional, because generically these disks have the critical points.

The following lemma from \cite{Kr6} provides that on extremal Teichm\"{u}ller  disks the functional $Z_n^0$ can be reconstructed from its matric $\ld_{Z_n^0}$.

\bigskip\noindent
{\bf Lemma 7}. {\it On any extremal Teichm\"{u}ller  disk
$\D(\psi_0) = \{t |\psi_0|/\psi_0: \ |t| < 1\} \subset \Belt(\D)_1$,
we have the equality
$$
\tanh^{-1}[\wt Z_n^0(F^{r|\psi_0|/\psi_0})] = \int\limits_0^r \ld_{\wt Z_n^0}(t) dt
$$
for each $r < 1$. }

\bigskip
In particular, this lemma implies the following result which is crucial for the proof of Zalcman's part of  Theorem 1.

\bigskip\noindent
{\bf Lemma 8}. \cite{Kr6} {\it For any $F(z) = z + b_0 + b_1 z^{-1} + \dots \in \Sigma_Q$,
we have the lower bound
 \be\label{19}
|\wt Z_n(F)| \ge \max_{b_0} |\wt Z_n(F_{b_0,b_1;1})|,
\end{equation}
where the maximum is taken over the set of admissible $b_0$. }

\bigskip\noindent
{\bf 3.4. Generalization}. The above lemmas, especially Lemma 7 on restoration of $Z_n$ by its infinitesimal form are also valid for functionals considered by Theorem 2, for example, $J(f) = a_n$ with $n > 2$.
This lemma deals with the extremal Teichm\"{u}ller disks in the sense of $\T$-equivalence, i.e.,
embedded in the universal space $\T$.

The second coefficient $a_2$ is somewhat specific. It is a part of the initial conditions
$w(0) = 0, \ w^\prime(0) = 1, \ w^{\prime\prime}(0) = 2 a_2$ defining the unique solution of the Schwarzian differential equation $S_w = \vp$ with given $\vp$ and intrinsically relates to quadratic differential $z^{-3} dz^2$ on $\hC$ having a simple pole at infinity. This differential expresses a
quasiconformal variation $a_2(w^\mu)$ via
$$
a_2(w^\mu) = - \fc{1}{\pi} \iint_{\D^*} \fc{\mu(\z)}{\z^3} d \xi d \eta + O(\|\mu\|_\iy^2),
\quad \|\mu\|_\iy \to 0.
$$
The indicated differential is holomorphic only on the punctured disk $\D^* \setminus \{\iy\}$, and the corresponding Teichm\"{u}ller disk $\{t z^3/|z|^3\}$ is extremal in the space $\T_1$, but not in $\T$.

The extremal dilatation along this disk is estimated by the hyperbolic distance of this punctured disk, which is greater that the distance on $\D^*$.

\bigskip\noindent
{\bf 3.5. Special quasiconformal deformations}. The following variational lemma is a special case of the general quasiconformal deformations constructed in \cite{Kr1}.

\bigskip\noindent
{\bf Lemma 9}. {\it Let $D$ be a simply connected domain on the Riemann sphere $\hC$. Assume that there are a set $E$ of positive two-dimensional Lebesgue measure and a finite number of points $z_1, z_2, ..., z_m$ distinguished in $D$. Let $\a_1, \a_2, ..., \a_m$ be non-negative integers assigned to $z_1,
z_2, ..., z_m$, respectively, so that $\a_j = 0$ if $z_j \in E$.

Then, for a sufficiently small $\ve_0 > 0$ and $\varepsilon \in (0, \varepsilon_0)$, and for any given collection of numbers $w_{sj}, s = 0, 1, ..., \a_j, \ j = 1,2, ..., m$ which satisfy the conditions
$w_{0j} \in D$, \
$$
|w_{0j} - z_j| \le \ve, \ \ |w_{1j} - 1| \le \ve, \ \ |w_{sj}| \le
\ve \ (s = 0, 1, \dots   a_j, \ j = 1, ..., m),
$$
there exists a quasiconformal automorphism $h$ of $D$ which is conformal on $D \setminus E$ and satisfies
$$
h^{(s)}(z_j) = w_{sj} \quad \text{for all} \ s =0, 1, ..., \a_j, \ j
= 1, ..., m.
$$
Moreover, the Beltrami coefficient $\mu_h(z) = \partial_{\ov z} h/\partial_z h$ of $h$ on $E$ satisfies $\| \mu_h \|_\iy \le M \ve$. The constants $\ve_0$ and $M$ depend only upon the sets $D, E$
and the vectors $(z_1, ..., z_m)$ and $(\a_1, ..., \a_m)$.

If the boundary $\partial D$ is Jordan or is $C^{l + \a}$-smooth, where $0 < \a < 1$ and $l \ge 1$, we can also take $z_j \in \partial D$ with $\a_j = 0$ or $\a_j \le l$, respectively.   }

\bigskip\bigskip
\centerline{\bf 4. PROOF OF THEOREM 1}

\bigskip
We accomplish the proof in five steps.

\bigskip\noindent
{\bf Step 1: Renormalization of functions}. First, we introduce the more general classes of univalent functions, which will be involved in the proof.

Taking into account Lemma 2, one can use the classes $S_{Q, \theta}(\D)$ of univalent functions in the disk $\D$ with expansions
$$
f(z) = e^{i \theta} z + a_2 z^2 + \dots, \quad - \pi \le \theta \le \pi,
$$
admitting quasiconformal extension to $\D^*$, and their subclasses $S_{z_0,\theta}$ consisting of $f \in S_{Q,\theta}$ with fix point at $z_0 \in \mathbb S^1$. The corresponding classes of univalent functions
$$
F(z) =  e^{- i \theta} z + b_0 + b_1 z^{-1} + b_2 z^{-2} + \dots .
$$
are denoted by $\Sigma_{Q, \theta}$ and $\Sigma_{z_0, \theta}$.
Consider their disjunct unions
$$
S^0 = \bigcup_{z_0 \in \mathbb S^1, \theta \in [-\pi,\pi]} S_{z_0,\theta}, \quad
\Sigma^0 = \bigcup_{z_0 \in \mathbb S^1, \theta \in [-\pi,\pi]} \Sigma_{z_0,\theta},
$$
and note that their closures $\ov{S^0}, \ \ov{\Sigma^0}$ in the topology of locally uniform convergence on the sphere $\C$ are compact.

Note that these families closely relates to the initial classes $S$ and $\Sigma$, because every $w \in S$ has its representative $\wh w$ in $\ov{S^0}$ (not necessarily unique) obtained by pre and post compositions of $w$ with rotations $z \mapsto e^{i \alpha} z$ about the origin, related by
$w_{\tau, \theta}(z) =  e^{- i \theta} w(e^{i \tau} z)$ with $\tau = \arg z_0$, where $z_0$ is a point for which $w(z_0) = e^{i \theta}$ is a common point of the unit circle and the  boundary of domain $w(\D)$. The existence of such a point follows from the classical Schwarz lemma.

Applying these classes instead of $S$ and $\Sigma$ involves a slight modification of the definition of Teichm\"{u}ller spaces $\T$ and $\T_1$ given above.
Actually, it yields the same model of the space $\T$.
As for the space $\T_1$, one must additionally add the third normalization condition that all
$F \in \Sigma_{z_0,\theta}$ preserve the point $z_0$ fixed. Thie provides a biholomorphically equivalent model of $\T_1$.

\bigskip\noindent
{\bf Step 2: Holomorphic lift of functionals $J_1(f)$ and $J_2(f)$ onto spaces $\T$ and $\T_1$}.
Using the relations between the coefficients $a_n$ of $f(z) \in S_{Q,\theta}$ and the corresponding coefficients $b_j$ of inversions $F_f(z) 1/f(1/z)$, given by
$$
b_0 + e^{2i \theta} a_2 = 0, \quad b_n + \sum \limits_{j=1}^{n}
\epsilon_{n,j}  b_{n-j} a_{j+1} + \epsilon_{n+2,0} a_{n+2} = 0,
\quad n = 1, 2, ... \ ,
$$
where $\epsilon_{n,j}$ are the entire powers of $e^{i \theta}$, one obtains successively the representations of $a_n$ by $b_j$:
 \be\label{20}
a_n = (- 1)^{n-1} \epsilon_{n-1,0}  b_0^{n-1} - (- 1)^{n-1} (n - 2)
\epsilon_{1,n-3} b_1 b_0^{n-3} + \text{lower terms with respect to} \ b_0.
\end{equation}
These relations transform the initial functionals $J_1(f^\mu), \ J_2(f^\mu)$ into the coefficient functionals $\wt J_1(F^\mu), \ J_2(F^\mu)$
on $\Sigma^0$ depending on the corresponding coefficients $b_j$. This dependence is holomorphic
from the Beltrami coefficients $\mu_F \in \Belt(\D)_1$ and from the Schwarzians $S_{F^\mu}$.

All functionals $J_1, J_2, J_{1,2}$ naturally extends to the generalized classes $S^0$ and $\Sigma^0$.

\bigskip
For any fixed $\theta$, the Taylor coefficients of functions $f^\mu \in S^0$ and $F^\mu \in \Sigma^0$ depend holomorphically on $\mu \in \Belt(\D)_1$ and on the Schwarzians $S_{F^\mu}$ as elements of $\B$. This
generates a holomorphic lifting the original functionals $J_s$ and $\wt J_s(F) \ (s = 1,2)$ onto the universal Teichm\"{u}ller space $\T \subset \B$ as holomorphic functions of $S_F \in \T$.

\bigskip
Our next goal is to lift these functionals onto the covering space $\T_1$. To reach this, we pass again to the functionals $\wh J_s(\mu) = \wt J_s(F^\mu)$. This relation lifts $J_s$ onto the ball $\Belt(\D)_1$.

Now we apply the $\T_1$-equivalence of maps $f^\mu$, i.e. the quotient map
$$
\phi_{\T_1}: \ \Belt(\D)_1 \to \T_1, \quad \mu \to [\mu]_{\T_1},
$$
which involves the homotopy of maps $F^\mu$ on the punctured disk $\D \setminus \{0\}$.
Thereby each functional $\wt J_s(F^\mu)$ is pushed down to a bounded holomorphic functional
$\mathcal J_s(X_{F^\mu})$
on the space $\T_1$, with the same range domain as $J_s$. Denote this functional by
$\mathcal J_s$.
\footnote{Recall that one has to deal with quadratic differentials $\vp(z) dz^2$ to insure the needed behaviour of the Schwarzian derivatives $\vp(h(z)) h^\prime(z)^2 = \vp(z)$
under conformal changes $h$ of variables.

This involves the equivalence classes of $\vp \in \B$ so that $\vp$ and $\vp_1$ are equivalent if $\vp_1(z) = \epsilon_1 \vp(\epsilon_2 z)$ for some constants $\epsilon_1, \ \epsilon_2$ with modulus $1$.
Such equivalence preserves $\B$ norm and moduli of coefficients.

The homogeneity of $J_s$ is compatible with this rotational invariance of $\T$.
}

Using the Bers isomorphism theorem, we regard the points of the space $\T_1 = \Fib (\T)$ as the pairs
$X_{F^\mu} = (S_{F^\mu}, F^\mu(0))$, where $\mu \in \Belt(\D)_1$ obey $\T_1$-equivalence.
So one obtains on the space $\Fib F(\T)$ the logarithmically plurisubharmonic functionals
$$
|\mathcal J_s(F^\mu, \ t)| = |\mathcal J_s(X_{F^\mu})|, \ s = 1, 2;
$$
and
\be\label{21}
J_{1,2}(F^\mu, t) = \max \left \{\max_S |J_1(F^\mu, t)|, \ \max_S |J_2((F^\mu,t)| \right\},
\quad t = F^\mu(0).
\end{equation}

\noindent
{\bf Step 3: Subharmonicity of maximal function generated by $J_{1,2}$}.
The functional (21) generates for any fixed $\theta \in [-\pi, \pi]$ and $F^\mu \in \Sigma_{Q,\theta}$ the maximal function
 \be\label{22}
u_\theta(t) = \sup_{S_{F^\mu}} |\mathcal J_{1,2}(S_{F^\mu}, t)|
\end{equation}
on the range domain $D_\theta$ of $F^\mu(0)$, taking the supremum over all $S_{F^\mu} \in \T$ admissible
for a given $t = F^\mu(0) \in D_\a$ (that means over the pairs $(S_{F^\mu}, t) \in \Fib(\T)$ with a fixed $t$).

The first crucial step in the proof of Theorem 1 is to establish that every function (22) inherits from $\mathcal J_1, \ \mathcal J_2$ subharmonicity in $t$, which we present as

\bigskip\noindent
{\bf Lemma 10}. {\it Every function $u_\theta(t)$ defiend by (22) with a fixed $\theta \in [-\pi, \pi]$ is logarithmically subharmonic in some domains $D_\theta$ located in the disk $\D_4 = \{|t| < 4\}$.
}

\bigskip\noindent
{\bf Proof}. Fix $\theta \in [-\pi, \pi]$ and, using the maps $F^\mu \in \Sigma_{Q,\theta}$,
apply a weak approximation of the underlying space $\T$ (and simultaneously of the space $\T_1$) by finite dimensional Teichm\"{u}ller spaces of the punctured spheres in the topology of locally uniform convergence
on $\C$.

Take the set of points
$$
E = \{e^{\pi s i/2^n}, \ s = 0, 1, \dots, 2^{n+1} - 1; \ n = 1, 2,
\dots\}
$$
(which is dense on the unit circle) and consider the punctured spheres
$$
X_m = \hC \setminus \{e^{\pi s i/2^n}, \ s = 0, 1, \dots, 2^{n+1} -
1\}, \quad m = 2^{n+1},
$$
and their universal holomorphic covering maps $g_m: \ \D \to X_m$ normalized by
$g_m(0) = 0, \ g_m^\prime(0) > 0$.

The radial slits from the infinite point to all the points $e^{\pi s i/2^n}$ form a canonical dissection $L_m$ of $X_m$ and define the simply connected surface $X_m^\prime = X_m \setminus L_m$. Any
covering map $g_m$ determines a Fuchsian group $\G_m$ of covering transformations uniformizing $X_m^\prime$, which act discontinuosly in both disks $\D$ and $\D^*$.

Every such group $G_m$ has a canonical (open) fundamental polygon $P_m$ of $\G_m$ in $\D$
corresponding to the dissection $L_m$. It is a regular circular $2^{n+1}$-gon centered at the origin
of the disk and can be chosen to have a vertex at the point $z = 1$.
The restriction of $g_m$ to $P_m$ is univalent, and as $m \to \iy$, these polygons entirely increase and exhaust the disk $\D$.

Similarly, we take in the complementary disk $\D^*$ the mirror polygons $P_m^*$ and the covering maps $g_m^*(z) = 1/\ov{g_m(1/ \ov z)}$ which define the mirror surfaces $X_m^*$.

Now we approximate the maps $F^\mu \in \Sigma_{Q, \theta}$ by homeomorphisms $F^{\mu_m}$ having in
$\D = \{|z| < 1\}$ the Beltrami coefficients
$$
\mu_m = [g_m]_{*} \mu := (\mu \circ g_m) \ov{g_m^\prime}/g_m^\prime, \ \ n = 1, 2, \dots \ .
$$
Each $F^{\mu_m}$ is again $k$-quasiconformal (where $k = \|\mu\|_\iy$)  and compatible with
the group $\G_m$. As $m \to \iy$, the coefficients $\mu_m$ are convergent to $\mu$ almost everywhere on $\C$; thus, the maps $F^{\mu_m}$ are convergent to $F^\mu$ uniformly in the spherical
metric on $\hC$.

Note also that $\mu_m$ depend holomorphically on $\mu$ as elements of $L_\iy$; hence, $F^{\mu_m}(0)$
is a holomorphic function of $t = F^\mu(0)$.

As a result, one obtains that the Beltrami coefficients
$$
\mu_{h,m} := [g_m]_{*} \mu_h
$$
and the corresponding values  $F^{\mu_{h,m}}(0)$ are holomorphic functions of the variable
$t = F^\mu(0)$.

By Hartogs theorem, each function $\mathcal J_s(S_{F^{\mu_m}}, t)$ with $t =
F^{\mu_m}(0)$ is jointly holomorphic in $(S_{F^{\mu_m}}, t) \in \mathcal F(\T)$.

We now choose in $\T(0, m) \setminus \{\mathbf 0\}$ represented as a subdomain of
the space $\B(\G_m)$ a countable dense subset
$$
E^{(m)} = \{\vp_1, \vp_2, \dots, \vp_p, \dots\}.
$$
For any of its point $\vp_p$,  the corresponding extremal
Teich\"{u}ller disk $\D(\vp_p)$ joining this point with the origin of $\B(\G_m)$ does
not meet other points from this set
(this follows from the uniqueness of Teich\"{u}ller extremal map). Recall also that
each disk $\D(\vp_p)$ is formed by the Schwarzians $S_{F^{\tau \mu_{p;m}}}$ with
$|\tau| < 1$ and
$$
\mu_{p;m}(z) = |\psi_{p;m}(z)|/\psi_{p;m}(z)
$$
with appropriate $\psi_{p;m} \in A_1(\D, \G_m), \ \|\psi_{p;m}\|_1 = 1$.

The restrictions of the functionals $\mathcal J_s(S_{F^{\tau \mu_{p;m}}}, t)$ to
these disks are holomorphic functions of $(\tau, t)$; moreover, the above construction
provides that all these restrictions are
holomorphic in $t$ in some common domain $D_m \subset \D_4$
containing the point $t = 0$, provided that $|\tau| \le k < 1$. We use the maximal
common holomorphy domain; it is located in a disk $\{|t| < r_0\}, \ r_0 < 4$.

Maximization over $\tau$ implies the logarithmically subharmonic functions
$$
U_{p;m}(t) = \sup_{|\tau| <1} |\mathcal J_{1,2}(S_{F^{\tau \mu_{p;m}}}, t)| \quad
(t = F^{\mu_{p;m}}(0), \ \ p = 1, 2, \dots)
$$
in the domain $D_m$. We consider the upper envelope of this sequence
$$
u_m(t) = \sup_p \ U_{p;m}(t)
$$
defined in some domain $D_m \subset \D_4$ containing the origin, and take its
upper semicontinuous regularization
$$
u_m(t) = \limsup\limits_{t^\prime \to t} u_m(t^\prime),
$$
which does not increase $\max |\mathcal J_{1,2}|$)
(by abuse of notation, we shall denote the regularizations by the same letter as the original functions).

Repeating this for all $m$, one obtains the sequences of monotone increasing functions $u_m(t)$ and of increasing domains $D_m$ exhausting a domain $D_\theta = \bigcup_m D_m$
such that each $u_m$ is subharmonic on $D_m$, and the limit function of this sequence is equal to the function (22). It is defined and subharmonic on the domain $D_\theta$. The lemma follows.

\bigskip\noindent
{\bf Step 4: Extremality of $\kp_\theta$ on $\Sigma_{\af}$}.
The upper semicontinuous envelope
$$
u(t) = \sup_\theta u_\theta(t)
$$
of functions (22) is logarithmically subharmonic in some domain $D_0$ containing the origin
$t = 0$, but gives only a rough upper bound
$$
\max_t u(t) > \max_S |J_{1,2}(f)|
$$
(in contrast to strongly rotationally homogeneous functionals). The weak rotational homogeneity
$|J_s(f_\a)| = |J_s(f)|$ for $f_\a(z) = e^{-i \a} f(e^{i \a} z)$ yields that the indicated domain
$D_0$ is a disk $\D_a = \{|t|< a\}$ of some radius $a \le 4$.

Our next goal is to find the best upper subharmonic dominant for $J_{1,2}$.
This will be done for the restriction of this functional to the set $\Sigma_{\af}$ by lifting to products of spaces $\T_1$.

\bigskip
First we establish the properties of the image of this set in the underlying space $\T$. Denote this image by $G_{\af}$.
Its structure is described by the following two lemmas. The first lemma reveals the intrinsic connection between the covering and distortion features of holomorphic functions.

\bigskip\noindent
{\bf Lemma 11}. \cite{Kr9} {\it Let $D$ be a bounded subdomain of $\C$,  $G$ be a domain in a complex Banach space $X = \{\mathbf x\}$ and $\chi$ be a holomorphic map from $G$ into the universal Teichm\"{u}ller space $\T = \Teich(D)$ with the base point $D$ (modeled as a bounded subdomain of $\B(D)$). Assume that $\chi(G)$ is a (pathwise connected) submanifold of finite or infinite dimension in $\T$.

Let $w(z)$ be a holomorphic univalent solution of the Schwarz differential equation
$$
S_w(z) = \chi(\x)
$$
on $D$ satisfying $w(0) = 0, \ w^\prime(0) = e^{i \theta}$ with the fixed $\theta \in [-\pi, \pi]$
and $\x \in G$ (hence $w(z) = e^{i \theta} z  + \sum_2^\infty a_n z^n$).
Put
\be\label{23}
|a_{2,\theta}^0| = \sup \{|a_2|: \ S_w \in \chi(G)\},
\end{equation}
and let $a_{2,\theta}^0 \ne 0$ and $w_0(z) = e^{i \theta} z + a_2^0 z^2 + \dots$ be one of the maximizing functions. Then:

(a) For every indicated function $w(z)$ , the image domain $w(D)$ covers entirely the disk
$D_{1/(2 |a_{2,\theta}^0|)} = \{|w| < 1/(2 |a_{2,\theta}^0|)\}$.

The radius value $1/(2 |a_{2, \theta}^0|)$ is sharp for this collection of
functions and fixed $\theta$, and the circle $\{|w| = 1/(2 |a_{2,\theta}^0|)$ contains points
not belonging to $w(\D)$ if and only if $|a_2| = |a_{2,\theta}^0|$
(i.e., when $w$ is one of the maximizing functions).

(b) The inverted functions
$$
W(\zeta) = 1/w(1/\zeta) = e^{i \theta}\zeta - a_2^0 + b_1 \zeta^{-1} + b_2 \zeta^{-2} + \dots
$$
with $\z \in D^{-1}$ map domain $D^{-1}$ onto a domain whose boundary is entirely
contained in the disk} $\{|W + a_{2,\theta}^0| \le |a_{2,\theta}^0|\}$.

\bigskip\noindent
{\bf Lemma 12}. {\it The (image of) set $\Sigma_{\af}$ is a three-dimensional subdomain in the space} $\T$.

\bigskip
This important lemma was already applied in \cite{Kr6}, \cite{Kr10}. Its proof is simple.
For any function $F_{b_0^0,b_1^0,t^0}$ with $|b_0^0| < 2, \ |b_1^0| < 1, \ |t^0| < 1$,
there is a small neighborhood $U_0$ of the point $(b_0^0,b_1^0,t^0)$ in $\C^3$ such that all functions
$F_{b_0,b_1,t}$ with $(b_0,b_1,t) \in U_0$ are univalent on $\D^*$ and admit quasiconformal extensions.
The map $t \mapsto S_{f_t}$is holomorphic in $\B$-norm, so the set $G_{\af}$ is
open in $\T$.

In addition, the line segment $[0, t^0]$ determines a curve in $G_{\af}$ joining the point
$F_{b_0^0,b_1^0,t^0}$ with the origin of $\T$; hence $G_{\af}$ is path-wise.

\bigskip
Denote the image of $\Sigma_{\af}$ in $\Fib(\T)$ by $\Fib (G_{\af})$; it is a complex four-dimensional subdomanifold
of $\Fib(\T)$. The restriction of functionals $\mathcal |J_s(S_{F^\mu}, t)|$ to this set are
plurisubharmonic on $\Fib(\Sigma_{\af})$, and their maximization similar to Step 3 implies a maximal
subharmonic function $u_{\af}(t)$ of type (22).

\bigskip
Our goal now is to find $\max_{\Sigma_{\af}} |\wt J_{1,2}(F)|$. Since the corresponding quantity
$$
a_{2,J} := \max_{F_f \in \Sigma_{\af}} |a_2(f)|
$$
is positive, one can apply Lemmas 11 and 12.

We select a dense subsequence $\{\theta_1, \theta_2, \dots\} \subset [-\pi, \pi]$ and define the corresponding functionals on the classes $S_{\theta_j}$ and $\Sigma_{\theta_j}$ replacing
the original functionals $\wt J_s$  as follows. Having the functions
$$
f_a(z) = e^{i \theta_m} z + a_2 z^2 + \dots
$$
and the corresponding
$$
F_a(z) = 1/f_a(z) = e^{- i \theta_m} z + b_0 + b_1 z^{-1} + \dots , \quad a =  e^{i \theta_m},
$$
consider $z^\prime = e^{i \theta_m} z$ as a new independent variable. Then
$$
f_a(z^\prime) = z^\prime + e^{- 2i \theta_m} (z^\prime)^2 + \dots, \quad
F_a(z^\prime) = z^\prime + b_0  e^{i \theta_m} + b_1  e^{i \theta_m} (z^\prime)^{-1} + \dots
$$
belong to $S_Q$ and $\Sigma_Q$ (in terms of variable $z^\prime$), and we set
$$
\wt J_{1,\theta_m}(F_a) = J_{1, \theta_m}(f_a) = \fc{e^{2i nm \theta} a_n^2 - e^{i(2n-1) m \theta} a_{2n-1}}{(n - 1)^2}, \quad
\wt J_{2,\theta_m}(F_a) = J_{2, \theta_m}(f_a) = \fc{e^{i nm \theta}a_n}{n}.
$$
This preserves the weak rotational homogeneity.

Note also that every $\wt J_{s,\theta_m}(F_a)$ obeys Lemma 8, replacing (19) by the inequality
$$
|\wt J_s(F_a)| \ge \max_{b_0} |\wt J_s(F_{b_0,b_1;1}^\theta)|,
$$
where
 \be\label{24}
F_{b_0,b_1;t}^\theta(z) = e^{i \theta} z + b_0 t + b_1 t^2 z^{-1} \quad \text{with} \ \ \theta = \arg a.
\end{equation}
The above relations result in
 \be\label{25}
\max_{\Sigma_{\theta_j}} |\wt J_{s,\theta_j}| = \max_\Sigma |\wt J_s| = \max_S |J_s|,
\end{equation}
and similarly for the corresponding collections $\Sigma_{\af,\theta_j}$ of functions (24).

\bigskip
Now consider the sequence of increasing products of the quotient spaces
 \be\label{26}
\mathcal T_m = \prod_{j=1}^m \ \wh \Sigma_{\theta_j}/\thicksim \
= \prod_{j=1}^m \{(S_{F_{\theta_j}}, F_{\theta_j}^{\mu_j}(0)) \} \ \simeq \T_1 \times \dots
\times \T_1,
\end{equation}
where the equivalence relation $\thicksim$ again means $\T_1$-equivalence.
The Beltrami coefficients  $\mu_j \in \Belt(\D)_1$ are chosen here independently. For any $\T_1$, presented in the right-hand side of (25), the corresponding values of $F_{\theta_j}^{\mu_j}(0)$ run over some domain $D_{\a_j} \subset \C$, and the corresponding collection
$\beta = (\beta_1, \dots, \beta_s)$
of the Bers isomorphisms
$$
\beta_j: \ \{(S_{W_{\theta_j}}, W_{\theta_j}^{\mu_j}(0))\} \to \mathcal F(\T)
$$
determines a holomorphic surjection of the space $\mathcal T_m$ onto the product of $m$ spaces
$\mathcal F(\T)$.

Letting
$$
\mathbf F_{\theta^\mu}(0) := (F_{\theta_1}^{\mu_1}(0), \dots , F_{\theta_m}^{\mu_m}(0)), \quad
\mathbf S_{\mathbf F_{\theta^\mu}} := (S_{F_{\theta_1}}, \dots, S_{F_{\theta_m}}),
$$
consider the holomorphic maps (vector-functions)
$$
\mathbf h(\mathbf S_{\mathbf F_{\theta}}) = (h_1(S_{F_{\theta_1}}), \dots h_m(S_{F_{\theta_m}})): \
\wh \Sigma_{\af, \theta} := \Sigma_{\af, \theta_1} \times \dots \times \Sigma_{\af, \theta_m} \to \C^m, \quad m = 1, 2, \dots,
$$
with
$$
h_j(S_{F_{\theta_j}}) = \wt Z_{n,\theta_j}(F_a), \quad j = 1, \dots, m,
$$
endowed with the polydisk norm
$$
\|\mathbf h\| = \max_j |h_j|
$$
on $\C^m$. Then by (25),
 \be\label{27}
\max_{\Sigma_{\af}} \|\mathbf h\| = \max_j \max_{\Sigma_{\af}}|h_j(S_{F_{\theta_j}})| = \max_S \  |J_{1,2}(f)|.
\end{equation}
The image of the set $\Sigma_{\af}$ under this embedding is the free product of $m$ factors $\Fib_j(G_{\af})$ with dimension $4m$. Note also that restriction of $\mathbf h$
to $\Sigma_{\af}$ is a polynomial map.

We now apply the construction from Step 3 simultaneously to each component $h(S_{F_{\theta_j}}, t)$
on the corresponding space $\T_1$ in (26) and obtain similar to Lemma 9 that the function
$$
u_m(t) = \max \bigl(|h(S_{F_{\theta_1}}, t)|, \dots, |h(S_{F_{\theta_m}}, t)|\bigr)
$$
is subharmonic in some domain $D_m$ containing the origin $t = 0$. This domain admits the rotational symmetry, hence it must be a disk $\D_{a_m}$ of some radius $a_m \le 4$.

This symmetry follows from rotational symmetry of the set $\Sigma_{\af}$, inherited by its images in
spaces $\T_1$ and $\mathcal T_m$, and from Lemma 9 (applied, if needed, to functions $F \in \Sigma_Q$ and a prescribed set $E$ in domain $F(\D)$) varying $F(0)$.

\bigskip
We apply this lemma to functions $F \in \Sigma_Q$ and take the prescribed set $E$ in domain $F(\D)$
to vary $F(0)$.

Each function $u_m(t)$ is a circularly symmetric function on its disk $\D_{a_m}$, and so is their upper envelope
$$
u_J(t) = \limsup_{m \to \iy} \wt u_m(t)
$$
(on some disk $\D_a, \ a > 0$). This envelope satisfies
$$
\max_{\D_a} u_J(t) = \max_S |J_{1,2}(f)|
$$
and  attains its maximal value at the boundary point $t = a$.

Noting that the closure of $\Sigma_{\af}$ contains the functions
$$
F_\theta(z) = z - 2 e^{i \theta} + e^{2i \theta} z^{-1}
$$
inverting the Koebe functions $\kp_\theta$, one derives that the radius $a$ must be equal $4$, which means that the range domain of $F^\mu(0)$ for $F^\mu \in \Sigma_{\af}$ coincides with the disk $\D_4$.
This implies that the  boundary points of this domain correspond only to functions $f(z)$ with $|a_2| = 2$, hence only to $\kp_\theta(z)$, and therefore, for all $f \in S$ with $F_f \in \Sigma_{\af}$, we have the estimate
 \be\label{28}
|J_{1,2}(f)| \le |J_{1,2}(\kp_\theta)| = 1.
\end{equation}
Moreover, the equality in the left hand part occurs only when $f = \kp_\theta$.

\bigskip\noindent
{\bf Step 5: Extremality of $\kp_\theta$ on the whole class $S$}.
It remains to establish that the relations (28) are also valid for all $f \in S$.

Noting that the homotopy disk of Koebe's function
$$
\D(\kp_\theta) = \{\kp_{\theta,t} = t^{-1} \kp_{\theta}(t z): \ |t| < 1\}
$$
is geodesic in universal Teichm\"{u}ller space and that the functionals $J_1, J_2$ generate on this
disk a conformal metric equal to the hyperbolic metric of the unit disk (which follows from Step 4),
one derives from asymptotic estimate (13) that this metric must be supporting at the origin for the
maximal infinitesimal metric $\ld_{J_{1,2}^0} (t)$ generated (via (19) applied separately to $J_1, J_2$)  along the extremal map $f_0$, since $\ld_{J_{1,2}^0} (0) \le \ld_\D(0) = 1$. Hence $f_0 = \kp_\theta$,
which proves the extremality of $\kp_\theta$ on the whole class $S$.

One also has that $\kp_\theta$ is unique extremal for $J_{1,2}$ on $S$, since  $\ld_{J_{1,2}^0} (t)$ majorates all conformal metrics determined by holomorphic maps $h: \D \to S$ on the corresponding disks $h(\D) \subset S$, and accordingly, the corresponding integral distances generated by metrics $\ld_{J_{1,2}^0}$ and $\ld_h$ must satisfy for any $f \in $S and any pair $(t, z)$ the inequality
$$
|J_{1,2}^0(f(t z))| \le |J_{1,2}^0(\kp_\theta(t z))| = (n - 1)^2 |t|^{2n-2}, \quad |t| \le 1,
$$
where the equality (even on one pair $(t, z)$) arises only when $f(z) = \kp_\theta(z)$.
This completes the proof of Theorem 1.

\bigskip
{\small\it{ \leftline{Department of Mathematics, Bar-Ilan
University, 5290002 Ramat-Gan, Israel}
\leftline{and Department of Mathematics, University of Virginia, Charlottesville, VA 22904-4137, USA}}}

\end{document}